\documentclass[12pt]{article}
\usepackage{amscd, amssymb, latexsym, amsmath}

\newtheorem{pro}{Proposition}
\newtheorem{lem}{Lemma}

\newenvironment{proof}
{\noindent {\em Proof}} {\hfill $\Box$}

\numberwithin{thm}{section} \numberwithin{cor}{section}
\numberwithin{pro}{section} \numberwithin{lem}{section}
\numberwithin{dfn}{section}
\numberwithin{rem}{section} \numberwithin{equation}{section}

\newcommand{\pyi}{\frac{\partial}{\partial y^i}}
\newcommand{\pyj}{\frac{\partial}{\partial y^j}}

\newcommand{\heat}{(\frac{d}{dt}-\Delta)}

\newcommand{\pyalpha}{\frac{\partial}{\partial y^\alpha}}
\newcommand{\pybeta}{\frac{\partial}{\partial y^\beta}}

\begin{document}
\title{Mean Curvature Flows and Isotopy of Maps Between Spheres}
\author{Mao-Pei Tsui \& Mu-Tao Wang }
%\address{Department of Mathematics\\
%Columbia University\\
%New York , NY 10027}
\date{February 19, 2003, revised September 4, 2003}
\maketitle
\centerline{email:\,tsui@math.columbia.edu,
mtwang@math.columbia.edu}

\begin{abstract}
Let $f$ be a smooth map between unit spheres of possibly different
dimensions.  We prove the global existence and convergence of the
mean curvature flow of the graph of $f$ under various conditions.
A corollary is that any area-decreasing map between unit spheres
(of possibly different dimensions) is isotopic to a constant map.
\end{abstract}

\section{Introduction}
Let $\Sigma_1$ and $\Sigma_2$ be two compact Riemannian manifolds
and $M = \Sigma_1 \times \Sigma_2$ be the product manifold. We
consider a smooth map $f:\Sigma_1
\rightarrow
\Sigma_2$ and denote the graph of $f$ by $\Sigma$; $\Sigma$ is a
submanifold of $M$ by the embedding $id\times f$.
 In \cite{wa1},
\cite{wa2}, and \cite{wa3}, the second author  studies the
deformation of $f$ by the mean curvature flow (see also the work
of Chen-Li-Tian \cite{clt}). The idea is to deform $\Sigma$ along
the direction of its mean curvature vector in $M$ with the hope
that $\Sigma$ will remain a graph.
%\marginlabel{Referee remark 2: Changed from the gradient flow}
  This is the negative gradient flow of the
volume functional and a stationary point is a ``minimal map"
introduced by Schoen in \cite{ sch}. In \cite{wa3}, the second
author  proves various long-time existence and convergence results
of graphical mean curvature flows in arbitrary codimensions under
assumptions on the Jacobian of the projection from $\Sigma$ to
$\Sigma_1$. This quantity is denoted by $*\Omega$ in \cite{wa3}
and $*\Omega>0$ if and only if $\Sigma$ is a graph over $\Sigma_1$
by the implicit function theorem. A crucial observation in
$\cite{wa3}$ is that $*\Omega$ is a monotone quantity under the
mean curvature flow when $*\Omega>\frac{1}{\sqrt{2}}$.

In this paper, we discover new positive geometric quantities
preserved by the graphical mean curvature flow. To describe these
results, we recall the differential of $f$, $df$, at each point of
$\Sigma_1$ is a linear map between the tangent spaces. The
Riemannian structures enables us to define the adjoint of $df$.
Let $\{\lambda_i\}$ denote the eigenvalues of $\sqrt{ (df)^T df}$,
or the singular values of $df$, where $(df)^T$ is the adjoint of
$df$. Note that $\lambda_i$ is always nonnegative. We say $f$ is
an {\it area decreasing map} if $\lambda_i \lambda_j < 1$ for any
$i\not=j$ at each point. In particular, $f$ is area-decreasing if
the $df$ has rank one everywhere. Under this condition, the second
author proves the Bernstein type theorem \cite{wa5} and interior
gradient estimates \cite{wa6} for solutions of the minimal surface
system. It is also proved in \cite{wa7} that the set of graphs of
area-decreasing linear transformations forms a convex subset of
the Grassmannian. We prove that this condition is preserved along
the mean curvature flow and the following global existence and
convergence theorem.

\vskip 10pt \noindent {\bf Theorem A.} {\it Let $\Sigma_1$ and
$\Sigma_2$ be compact Riemannian manifolds of constant curvature
$k_1$ and $k_2$ respectively. Suppose $k_1\geq |k_2|$, $k_1+k_2
>0$ and $dim(\Sigma_1) \ge 2$. If $f$ is a smooth area decreasing map from $\Sigma_1$ to $\Sigma_2$, the mean curvature flow
of the graph of $f$
 remains the graph of an area decreasing map, exists for all time, and converges smoothly to the graph of
  a constant map.} \vskip 10pt

We remark that the condition $k_1\geq |k_2|$ is enough to prove
the long time existence of the flow. The following is an
application to determine when a map between spheres is
homotopically trivial. \vskip 10pt \noindent {\bf Corollary A }
{\it Any area-decreasing map from $S^n$ to $S^m$ with $n\geq 2$ is
homotopically trivial.}
 \vskip 10pt
When $m=1$, the area-decreasing condition always holds and the
above statement follows from the fact that $\pi_n(S^1)$ is trivial
for $n\ge2$.
 We remark that
the result when $m=2$ is proved by the second author in \cite{wa4}
using a somewhat different method. The higher homotopy groups
$\pi_n(S^m)$ has been computed in many cases and it is known that
homotopically nontrivial maps do exist when $n\geq m$. Since an
area-decreasing map may still be surjective when $n>m$, we do not
know any topological method that would imply such a conclusion.

We would like to thank Professor R. Hamilton, Professor D. H.
Phong and Professor S.-T. Yau for their constant advice,
encouragement and support.

\section{Preliminaries}
In this section, we recall notations  and formulae for mean
curvature flows. Let $f:\Sigma_1\rightarrow \Sigma_2$ be a smooth
map between Riemannian manifolds. The graph of $f$ is an embedded
submanifold $\Sigma$ in $M=\Sigma_1\times \Sigma_2$.  At any point
of $\Sigma$, the tangent space of $M$, $TM$ splits into the direct
sum of the tangent space of $\Sigma$, $T\Sigma$ and the normal
space $N\Sigma$, the orthogonal complement of the tangent space
$T\Sigma$ in $TM$. There are
 isomorphisms $T\Sigma_1\rightarrow T\Sigma$
by $X\mapsto X+df(X)$ and $T\Sigma_2\rightarrow N\Sigma$ by $Y
\mapsto Y-(df)^T(Y)$ where $(df)^T: T\Sigma_2 \rightarrow
T\Sigma_1$ is the adjoint of $df$.

We assume the mean curvature flow of $\Sigma$ can be written as a
graph of $f_t$ for $t\in [0, \epsilon)$ and derive the equation
satisfied by $f_t$. The mean curvature flow is given by a smooth
family of immersions $F_t$ of $\Sigma$ into $M$ which satisfies

\[(\frac{\partial F}{\partial t})^\perp
=H\] where $H$ is the mean curvature vector in $M$ and
$(\cdot)^\perp$ denotes the projection onto the normal space
$N\Sigma$. Notice that we do not require $\frac{\partial
F}{\partial t}$ is in the normal direction since the difference is
only a tangential diffeomorphism (see for example White \cite{wh}
for the issue of parametrization). By the definition of the mean
curvature vector, this equation is equivalent to
\[(\frac{\partial F}{\partial t})^\perp
=(\Lambda^{ij} \nabla_{\frac{\partial F} {\partial
x^i}}^M\frac{\partial F}{\partial x^j})^ \perp\] where
$\Lambda^{ij}$ is the inverse to
 the induced metric $\Lambda_{ij}=\langle \frac{\partial F}{\partial x^i},
 \frac{\partial F}{\partial x^j}\rangle$
on $\Sigma$.

In terms of coordinates $\{y^A\}_{A=1\cdots n+m}$ on $M$, we have
 \[\Lambda^{ij} \nabla^M_{\frac{\partial F}{\partial
x^j}} {\frac{\partial F}{\partial
x^i}}=\Lambda^{ij}(\frac{\partial^2 F^A}{\partial x^i\partial
x^j}+\frac{\partial F^B}{\partial x^i} \frac{\partial
F^C}{\partial x^j} \Gamma_{BC}^A) \frac{\partial}{\partial y^A}\]
where $\Gamma_{BC}^A$ is the Christoffel symbol of $M$ and thus

\[(\Lambda^{ij} \nabla^M_{\frac{\partial F}{\partial
x^j}} {\frac{\partial F}{\partial
x^i}})^\perp=\Lambda^{ij}(\frac{\partial^2 F^A}{\partial
x^i\partial x^j}+\frac{\partial F^B}{\partial x^i} \frac{\partial
F^C}{\partial x^j}
\Gamma_{BC}^A-\tilde{\Gamma}_{ij}^k\frac{\partial F^A}{\partial
x^k}) \frac{\partial}{\partial y^A}
\]where $\tilde{\Gamma}_{ij}^k$ is the Christoffel symbol of the induced metric on $\Sigma$.

 By assumption, the embedding is given by the graph of $f_t$.
We fix a coordinate system $\{x^i\}$ on $\Sigma_1$ and consider
$F:\Sigma_1\times [0, T)\rightarrow M$ given by

\[F(x^1, \cdots, x^n, t)=(x^1,\cdots,x^n, f^{n+1}, \cdots, f^{n+m}).\]
We shall use $i,j, k, l\cdots=1\cdots n$ and $\alpha, \beta,
\gamma=n+1\cdots n+m$ for the indices. Of course $f^\alpha
=f^\alpha(x^1,\cdots, x^n,t)$ is time-dependent.

Therefore $\frac{\partial F}{\partial t}=\frac{\partial
f^\alpha}{\partial t}\pyalpha$ and

\[\Lambda^{ij}(\frac{\partial^2 F^A}{\partial x^i\partial x^j}
+\frac{\partial F^B}{\partial x^i} \frac{\partial F^C}{\partial
x^j} \Gamma_{BC}^A )\frac{\partial}{\partial
y^A}=\Lambda^{ij}(\frac{\partial^2 f^\alpha}{\partial x^i\partial
x^j} \frac{\partial}{\partial
y^\alpha}+\Gamma_{ij}^l\frac{\partial}{\partial
y^l}+\frac{\partial f^\beta}{\partial x^i} \frac{\partial
f^\gamma}{\partial x^j} \Gamma_{\beta\gamma}^\alpha
\frac{\partial}{\partial y^\alpha}).\]

Thus the mean curvature flow equation is equivalent to the normal
part of

\[[ \frac{\partial f^\alpha}{\partial t}-\Lambda^{ij}(\frac{\partial^2
f^\alpha}{\partial x^i\partial x^j} +\frac{\partial
f^\beta}{\partial x^i} \frac{\partial f^\gamma}{\partial x^j}
\Gamma_{\beta\gamma}^\alpha)]\frac{\partial}{\partial y^\alpha}-
\Lambda^{ij}\Gamma_{ij}^l\frac{\partial}{\partial y^l}\] is zero.

Now given any vector $a^i\pyi+b^\alpha \pyalpha$, the equation
that the normal part being zero is equivalent to

\begin{equation}b^\alpha -a^i \frac{\partial f^\alpha}{\partial
x^i}=0\end{equation} for each $\alpha$.  Therefore we obtain the
evolution equation for $f$

\begin{equation}\label{equation_for_f}\frac{\partial f^\alpha}{\partial
t}-\Lambda^{ij}(\frac{\partial^2 f^\alpha}{\partial x^i\partial
x^j} +\frac{\partial f^\beta}{\partial x^i} \frac{\partial
f^\gamma}{\partial x^j}
\Gamma_{\beta\gamma}^\alpha+\Gamma_{ij}^k\frac{\partial f^\alpha
}{\partial x^k})=0.\end{equation} where $\Lambda^{ij}$ is the
inverse to $g_{ij}+h_{\alpha\beta}\frac{\partial
f^\alpha}{\partial x^i}\frac{\partial f^\beta}{\partial x^j}$ and
$g_{ij}=\langle \pyi, \pyj \rangle$ and
$h_{\alpha\beta}=\langle\pyalpha , \pybeta\rangle$  are the
Riemannian metrics on $\Sigma_1$ and $\Sigma_2$, respectively.
$\Gamma_{ij}^k$ and $\Gamma_{\beta\gamma}^\alpha$ are the
Christoffel symbols of $g_{ij}$ and $h_{\alpha\beta}$
respectively.

(\ref{equation_for_f}) is a nonlinear parabolic system and the
usual derivative estimates do not apply to this equations.
However, the second author  in \cite{wa3} identifies a geometric
quantity in terms of the derivatives of $f^\alpha$ that satisfies
the maximum principle; this quantity and its evolution equation
are recalled in the next section.

\section{Two evolution equations}

In this section, we recall two evolution equations along the mean
curvature flow. The basic set-up is a mean curvature flow
$F:\Sigma\times[0,T)\rightarrow  M$ of an $n$ dimensional
submanifold $\Sigma$ inside an $n+m$ dimensional Riemannian
manifold $M$. Given any parallel tensor on $M$, we may consider
the pull-back tensor by $F_t$ and consider the evolution equation
with respect to the time-dependent induced metric on
$F_t(\Sigma)=\Sigma_t$. For the purpose of applying maximum
principle, it suffices to derive the equation at a space-time
point. We write all geometric quantities in terms of orthonormal
frames keeping in mind all quantities are defined independent of
choices of frames. At any point $p\in \Sigma_t$, we choose any
orthonormal frames $\{e_i\}_{i=1\cdots n}$ for $T_p\Sigma_t$ and
$\{e_{\alpha}\}_{\alpha=n+1 \cdots n+m}$ for $N_p\Sigma_t$. The
second fundamental form $h_{\alpha ij}$ is denoted by $h_{\alpha
ij} =\langle\nabla_{e_i}^M e_j, e_{\alpha}\rangle$ and the mean
curvature vector is denoted by $H_\alpha=\sum_ih_{\alpha ii}$. For
any $j,k$, we pretend
\[h_{n+i, jk}=0\] if $i>m$.

 When $M=\Sigma_1\times \Sigma_2$ is the product of $\Sigma_1$ and $\Sigma_2$, we denote the projections by
$\pi_1:M\rightarrow \Sigma_1$ and $\pi_2:M\rightarrow \Sigma_2$.
By abusing notations, we also denote the differentials by
$\pi_1:T_p M\rightarrow T_{\pi_1(p)}\Sigma_1$ and $\pi_1:T_p
M\rightarrow T_{\pi_2(p)}\Sigma_2$ at any  point $p\in M$. The
volume form $\Omega$ of $\Sigma_1$ can be extended to a parallel
$n$-form on $M$. For an oriented orthonormal basis $e_1 \cdots
e_n$ of $T_p \Sigma$, $\Omega(e_1,\cdots,
e_n)=\Omega(\pi_1(e_1),\cdots ,\pi_1(e_n))$ is the Jacobian of the
projection from $T_p\Sigma$ to $T_{\pi_1(p)}\Sigma_1$. This can
also be considered as the pairing between the $n$-form $\Omega$
and the $n$-vector $e_1\wedge\cdots\wedge e_n$ representing
$T_p\Sigma$. We use $*\Omega$ to denote this function as $p$
varies along $\Sigma$. By the implicit function theorem,
$*\Omega>0$ at $p$ if and only if $\Sigma$ is locally a graph over
$\Sigma_1$ at $p$. The evolution equation of $*\Omega$ is
calculated in Proposition 3.2 of \cite{wa3}.

When $\Sigma$ is the graph of $f:\Sigma_1\rightarrow \Sigma_2$,
the equation at each point can be written in terms of singular
values of $df$ and special bases adapted to $df$. Denote the
singular values of $df$, or eigenvalues of $(df)^T df$, by
$\{\lambda_i\}_{i=1\cdots n}$. Let $r$ denote the rank of $df$.
 We can rearrange them so that
$\lambda_i=0$ when $i$ is greater than $r$. By singular value
decomposition, there exist orthonormal bases $\{a_i\}_{i=1 \cdots
n}$ for $T_{\pi_1(p)}\Sigma_1$ and $\{a_{\alpha}\}_{
\alpha=n+1\cdots n+m}$ for $T_{\pi_2(p)} \Sigma_2$ such that
\[df(a_i)=\lambda_i a_{n+i}\,\,\] for $i$ less than or equal to $r$
and $df(a_i)=0$ for $i$ greater than $r $. Moreover,
\begin{equation}\label{onf1_1} e_i=
\begin{cases} \frac{1}{\sqrt{1+\lambda_i^2}}(a_i+\lambda_i
a_{n+i})& \text{if}\,\, {1 \le i \le r   }\\
 a_i & \text{if}\,\, {r+1 \le i \le n   }
 \end{cases}
\end{equation}
 becomes an orthonormal basis for $T_p\Sigma$
 and
 \begin{equation}\label{onf1_3}
e_{n+p} =
\begin{cases}
\frac{1}{\sqrt{1+\lambda_p^2}}(a_{n+p} -\lambda_p
a_{p})& \text{if}\,\,1 \le p \le r\\
a_ {n+p}& \text{if}\,\,r+1 \le p \le m
\end{cases}
\end{equation}
becomes an orthonormal basis for $N_p\Sigma$.

In terms of the singular values $\lambda_i$,
\begin{equation}\label{Omega}*\Omega=\frac{1}{\sqrt{\prod_{i=1}^n(1+\lambda_i^2)}}\end{equation}
With all the notations understood, the following result is
essentially derived in Proposition 3.2 of \cite{wa3}
 by noting that $(\ln *\Omega)_k = -(\sum_i \lambda_i h_{n+i, i k})$.

\begin{pro}\label{equation}
Suppose $M=\Sigma_1\times \Sigma_2$ and $\Sigma_1$ and $\Sigma_2$
are compact Riemannian manifolds of constant curvature $k_1$ and
$k_2$ respectively. With respect to the particular bases given by
the singular value decomposition of $df$, $\ln *\Omega$ satisfies
the following equation.
\begin{equation}\label{lneq}
\begin{split}
&\heat\ln *\Omega =\sum_{\alpha, i, k} h_{\alpha ik}^2 + \sum_{k,
i} \lambda_{ i}^2 h_{n+i,ik}^2
+2\sum_{k, i<j}\lambda_{i}\lambda_{j} h_{n+j, ik} h_{n+i, jk}\\
&+\sum_{ i}\frac{\lambda_{ i}^2}{1+ \lambda_i^2}\Big[(k_1+k_2)
(\sum_{j\not= i}\frac{1}{1+\lambda_j^2}) +k_2 (1-n)\Big]
\end{split}
\end{equation}
\end{pro}

Next we recall the evolution equation of parallel two tensors from
\cite{sw}. The calculation indeed already appears in \cite{wa1}.
The equation will be used later to obtain more refined
information.
 Given a parallel two-tensor $S$ on
$M$, we consider the evolution of $S$ restricted to $\Sigma_t$.
This is a family of time-dependent symmetric two tensors on
$\Sigma_t$.

\begin{pro}\label{evo_sym_l}
Let $S$ be a parallel two-tensor on $M$.   Then the pull-back of
$S$ to $\Sigma_t$ satisfies the following equation.
\begin{equation}\label{evo_symb}
\begin{split}
\heat S_{ij}&= - h_{\alpha i l}H_{\alpha}S_{lj}
- h_{\alpha j l}H_{\alpha}S_{li} \\
& + R_{k i k \alpha}S_{\alpha j} +  R_{k j k\alpha }S_{\alpha i}\\
& +h_{\alpha k l}h_{\alpha k i}S_{lj}
  +h_{\alpha k l}h_{\alpha k j}S_{li}
-2h_{\alpha k i} h_{\beta k j}S_{\alpha \beta}
\end{split}
\end{equation}
 where $\Delta $
is the rough Laplacian on two-tensors over $\Sigma_t$ and
$S_{\alpha i}=S(e_\alpha, e_i)$, $S_{\alpha \beta}=S(e_{\alpha},
e_{\beta})$, and $R_{k i k \alpha}=R(e_k, e_i, e_k, e_\alpha)$ is
the curvature of $M$.
\end{pro}

The  evolution equations (\ref{evo_symb}) of $S$ can be written in
terms of  evolving orthonormal frames as in Hamilton \cite{ha3}.
If the orthonormal frames
\begin{equation}\label{onf}
F=\{F_1, \cdots, F_a, \cdots,F_{n}\}
\end{equation}
are given in local coordinates by
\[
F_a = F_a^i \frac{\partial }{\partial x_i}\; .
\]
To keep them orthonormal, i.e. $ g_{ij}F_a^iF_b^j = \delta_{ab}$,
we evolve $F$ by the formula
\[
\frac{\partial }{\partial t} F_a^i = g^{ij}g^{\alpha\beta}
h_{\alpha j l}H_{\beta}F_a^l \, .
\]

Let $S_{ab} = S_{ij}F_a^iF_b^j$ be the components of $S$ in $F$.
Then $S_{ab}$ satisfies the following equation
\begin{equation}\label{onframes}
\begin{split}
\heat S_{ab}&=  \
 R_{c a c \alpha}S_{\alpha b} +  R_{c b c \alpha}S_{\alpha a}\\
& \ +h_{\alpha c d}h_{\alpha c a}S_{db}
  +h_{\alpha c d}h_{\alpha c b}S_{da} \\
& \ -2h_{\alpha c a} h_{\beta c b}S_{\alpha \beta} \, .
\end{split}
\end{equation}

\section{Preserving the distance-decreasing condition}
In this section, we show the condition $|df|<1$, or each singular
value $\lambda_i<1$, is preserved by the mean curvature flow. This
result will not be used in proof of the Theorem A. But the proof
of  Theorem A depends on the computation in this section. The
tangent space of $M=\Sigma_1 \times \Sigma_2$ is identified with
$T\Sigma_1 \oplus T\Sigma_2$. Let $\pi_1$ and $\pi_2$ denote the
projection onto the first and second summand in the splitting. We
define the parallel symmetric two-tensor $S$  by
\begin{equation}\label{sym}
S(X,Y)=\langle \pi_1(X), \pi_1(Y)\rangle\ -\langle \pi_2(X),
\pi_2(Y)\rangle
\end{equation}
 for any
 $X, Y \in TM $.

Let $\Sigma$ be the graph of $f:\Sigma_1\rightarrow \Sigma_1
\times \Sigma_2$.  $S$ restricts to a symmetric two-tensor on
$\Sigma$ and we can represent $S$ in terms of the orthonormal
basis (\ref{onf1_1}).

Let $r$ denote the rank of $df$.
By (\ref{onf1_1}), it is not hard to check
\begin{equation}\label{pi1}
\begin{split}
\pi_1(e_i) & = \ \frac{a_i} {\sqrt{1+{\lambda_i}^2}} \; ,
 \pi_2(e_i)= \ \frac{\lambda_ia_{n+i}} {\sqrt{1+{\lambda_i}^2}}
\,\, \
\text{for}\,\,  1 \le i \le r \,\, , \\
\text{and} \ \pi_1(e_i) & = \ a_i \;,
 \pi_2(e_i)= \ 0
\,\,
\text{for}\,\,  r+1 \le i \le n.
\end{split}
\end{equation}
Similarly, by (\ref{onf1_3}) we have
 \begin{equation}\label{pi1_1}
\begin{split}
\pi_1(e_{n+p}) &= \ \frac{-\lambda_p a_p} {\sqrt{1+{\lambda_p}^2}}\; ,
\pi_2(e_{n+p}) = \ \frac{a_{n+p}}
{\sqrt{1+{\lambda_p}^2}}\,\,
\text{for}\,\,  1 \le p \le r\;,\,\, \\
\,\,
\text{and}\,\,
\pi_1(e_{n+p}) &= 0 \; ,
\pi_2(e_{n+p}) = a_{n+p}  \,\,
\text{for}\,\,r+1 \le p \le m\; .
\end{split}
\end{equation}
From the definition of $S$, we have
\begin{equation}\label{S_ij}
\begin{split}
 S(e_i, e_j)
 =  \frac{1- \lambda_i^2}{1+{\lambda_i}^2} \delta_{ij}\; .
\end{split}
\end{equation}
In particular, the eigenvalues of $S$ are
\begin{equation}\label{eigenvalues}
 \frac{1-{\lambda_i}^2 }{1+{\lambda_i}^2}, \,i= 1 \cdots n\; .
\end{equation}
 Notice that $S$  is positive-definite if and only if
\[
 |\lambda_i| < 1
\] for any singular value $\lambda_i$ of $df$.

Now,  at each point we express $S$ in terms of the orthonormal
basis $\{{e_i}\}_{i = 1 \cdots n}$ and  $\{{e_{\alpha}}\}_{\alpha
= n+1
\cdots n+m}$. Let $I_{k \times k}$ denote a $k$ by $k$ identity
matrix. Then $S$ can be written in the block form
\begin{equation}\label{block}
 S = \Big(
S(e_k,e_l)
\Big)_{ 1 \le k,l \le {n+m}}=
 \left(
  \begin{matrix}
 B & 0 & D & 0 \\
0 & I_{n-r \times n-r} & 0 & 0 \\
D & 0 & -B & 0 \\
0 & 0 &0 &  - I_{m-r \times m-r}
\end{matrix}
\right)
\end{equation}
where $B$ and $D$ are  $r$ by $r$ matrices with $ B_{ij}=
S(e_i,e_j) = \frac{1-\lambda_i^2}{1+\lambda_i^2}\delta_{ij}$ and $
D_ {ij}= S(e_i,e_{n+j}) = \frac{-2\lambda_i}{1+\lambda_i^2}
\delta_{ij} $ for  $1 \le i, j \le r$. We show that the positivity
of $S$ is preserved by the mean curvature flow. We remark that a
similar positive definite tensor has been considered for the
Lagrangian mean curvature flow in Smoczyk \cite{sm} and
Smoczyk-Wang \cite{sw}. The following lemma shows that the
distance decreasing condition is preserved by the mean curvature
flow if $k_1\geq |k_2|$.

\begin{lem}
The condition
\begin{equation}\label{covexityeq}
T_{ij} = S_{ij} - \epsilon g_{ij} > 0 \, \, \text{ for some
}\epsilon \ge 0
\end{equation}
 is preserved by the mean
curvature flow if $k_1\geq |k_2|$.
\end{lem}

\begin{proof}.
We compute the evolution equation for $T_{ij}$. From Proposition
(\ref{evo_sym_l}) and
\[
\frac{\partial}{\partial t}g_{ij} = -2  h_{\alpha i j}H_{\alpha}
\, ,
\]
 we have
\begin{equation}\label{evol_nor2}
\begin{split}
\heat T_{ij}&= - h_{\alpha i l}H_{\alpha}T_{lj} - h_{\alpha j
l}H_{\alpha}T_{li}
 + R_{k i k \alpha}S_{\alpha j} +  R_{k j k\alpha }S_{\alpha i}\\
& +h_{\alpha k l}h_{\alpha k i}T_{lj}
  +h_{\alpha k l}h_{\alpha k j}T_{li}+ 2\epsilon h_{\alpha k i}h_{\alpha k j}
-2h_{\alpha k i} h_{\beta k j}S_{\alpha \beta}.
\end{split}
\end{equation}

To apply Hamilton's maximum principle, it suffices to prove that
$N_{ij} V^i V^j \ge 0$ for any null eigenvector $V$ of $T_{ij}$,
where $ N_{ij} $ is the right hand side of (\ref{evol_nor2}).
Since $V$ is a null eigenvector of $T_{ij}$, it satisfies $\sum_j
T_{ij}V^j = 0$ for any $i$, and thus $ N_{ij} V^i V^j$ is equal to
\begin{equation}\label{null}\\
\begin{split}
2 \epsilon h_{\alpha k i}h_{\alpha k j} V^i V^j +
 2 R_{k i k \alpha}S_{\alpha j} V^i V^j
-2h_{\alpha k i} h_{\beta k j}S_{\alpha \beta}V^i V^j\; .
\end{split}
\end{equation}

Obviously, the first term of (\ref{null}) is nonnegative. Applying
the relation in (\ref{block})  to the last term of (\ref{null}) we
obtain
\[
-2h_{\alpha k i} h_{\beta k j}S_{\alpha \beta}V^i V^j =
\sum_{1 \le p, q \le  r}2h_{n+p k i} h_{n+q k j}S_{pq}V^i V^j +
\sum_{r+1 \le p, q \le  m}2h_{n+p k i} h_{n+q k j}V^i V^j\, .
\]
 Since
$T_{pq} \ge 0$ implies that $S_{pq} \ge \epsilon g_{pq}$, we
obtain $-2h_{\alpha k i} h_{\beta k j}S_{\alpha \beta}V^i V^j \ge
0$. In the next lemma we  show that $ R_{kik \alpha} S_{ \alpha
j}$ is nonnegative definite whenever $ S_{ij}$ is under the
curvature assumption $k_1\ge |k_2|$.

\end{proof}

\begin{lem}
\begin{equation}\label{curv_vec}
\begin{split}
  R_{kik \alpha} S_{\alpha
j}=\frac{\lambda_i^2}{(1+\lambda_i^2)^2}\left[
{(k_1-k_2)(n-1)}+(k_1+k_2)
 \sum_{k\not= i}\frac{1-\lambda_k^2}{1+\lambda_k^2} \right]\delta_{ij}.
\end{split}
\end{equation}
\end{lem}

\begin{proof}. We follow the calculation of the curvature terms in
$\cite{wa3}$.

\[
\begin{split}
& \quad \sum_k R(e_{\alpha}, e_k, e_k ,e_i)\\
& =  \sum_k  R_1(\pi_1(e_{\alpha}), \pi_1(e_k), \pi_1(e_k)
,\pi_1(e_i)) +R_2(\pi_2(e_{\alpha}), \pi_2(e_k), \pi_2(e_k)
,\pi_2(e_i))\\
& =  \sum_k  k_1\Big[\langle\pi_1(e_{\alpha}), \pi_1(e_k)\rangle
\langle\pi_1(e_k), \pi_1(e_i)\rangle -\langle\pi_1(e_{\alpha}),
\pi_1(e_i)\rangle
\langle\pi_1(e_k), \pi_1(e_k)\rangle\Big]\\
& \quad + k_2\Big[\langle\pi_2(e_{\alpha}), \pi_2(e_k)\rangle
\langle\pi_2(e_k), \pi_2(e_i)\rangle -\langle\pi_2(e_{\alpha}),
\pi_2(e_i)\rangle \langle\pi_2(e_k), \pi_2(e_k)\rangle\Big] \; .
\end{split}
\]

Notice that $\langle\pi_2(X), \pi_2(Y)\rangle=\langle
X,Y\rangle-\langle\pi_1(X), \pi_1(Y)\rangle $ since $T\Sigma_1
\perp T\Sigma_2$.
Therefore
\[
\begin{split}
& \quad \sum_k R(e_{\alpha}, e_k, e_k ,e_i)\\
& = \sum_k (k_1+k_2)\Big[\langle\pi_1(e_{\alpha}),
\pi_1(e_k)\rangle \langle\pi_1(e_k), \pi_1(e_i)\rangle
-\langle\pi_1(e_{\alpha}), \pi_1(e_i)\rangle |\pi_1(e_k)|^2\Big]\\
& \quad +k_2(n-1)\langle\pi_1(e_{\alpha}), \pi_1(e_i)\rangle\\
\end{split}
\]

Now use $\pi_1(e_\alpha)= -\lambda_p \pi_1(e_p) \delta_{\alpha,
n+p}$ and  $S(e_{j}, e_{n+p})
 =  -\tfrac{ 2\lambda_j
 \delta_{jp}}{1+{\lambda_j}^2}$ in (\ref{block}),
 we have
\[
\begin{split}
& \ \sum_{\alpha,k} R_{kik\alpha} S_{\alpha j}
= -\sum_{p,k} R_{n+p,kki} S_{n+p, j} \\
= & \sum_{p,k}
\Big\{\lambda_{p}(k_1+k_2)\left[\langle\pi_1(e_{p}),
\pi_1(e_k)\rangle \langle\pi_1(e_k), \pi_1(e_i)\rangle
-\langle\pi_1(e_{p}), \pi_1(e_i)\rangle |\pi_1(e_k)|^2\right]\\
& \quad +\lambda_{p}k_2(n-1)\langle\pi_1(e_{p}), \pi_1(e_i)
\rangle\Big\} S_{n+p, j} \\
= & -\frac{ 2\lambda_i^2 }{1+{\lambda_i}^2}
\Big\{(k_1+k_2)\Big[\frac{\delta_{ij}
}{(1+{\lambda_i}^2)^2}
-\frac{\delta_{ij}  }{1+\lambda_i^2}\sum_{k}|\pi_1(e_k)|^2\Big]\\
& \quad +k_2(n-1)\frac{\delta_{ij} }{1+{\lambda_i}^2} \Big\}.
\end{split}
\]
Recall that $ |\pi_1(e_k)|^2=\frac{1 }{1+\lambda_k^2}$ and we
obtain
\[
\begin{split}
 R_{kik \alpha} S_{\alpha j}
=  \frac{2\lambda_{i}^2\delta_{ij}}{(1+\lambda_{i}^2)^2}
\left[(k_1+k_2) (\sum_{k \neq i}\frac{1}{1+\lambda_{k}^2}) +
     k_2 (1-n )\right].\\
\end{split}
\]
This can be further simplified by noting
\begin{equation}\label{useful_identity}
\begin{split}
& (k_1+k_2) (\sum_{k \neq i}\frac{1}{1+\lambda_{k}^2}) +
     k_2 (1-n )={\frac{(k_1-k_2)(n-1)}{2}}+(k_1+k_2)
 \sum_{k\not= i}\frac{1-\lambda_k^2}{2(1+\lambda_k^2)}
\end{split}
\end{equation}
where we  use the following identity for each $i$
\[
\begin{split}
& (\sum_{k\not= i}\frac{1}{1+\lambda_k^2})- \frac{n-1}{2} =
\sum_{k\not= i}(\frac{1}{1+\lambda_k^2}-\frac{1}{2})=\sum_{k\not=
i}\frac{1-\lambda_k^2}{2(1+\lambda_k^2)}.
\end{split}
\]
\end{proof}

\section{Preserving the area-decreasing condition}
In this section, we  show that the area decreasing condition is
preserved along the mean curvature flow. In the following, we
 require that $n=dim(\Sigma_1) \ge 2$. By (\ref{eigenvalues}),
the sum of any two eigenvalues of $S$ is
\begin{equation}\label{eigs}
 \frac{1- \lambda_i^2}{1+ \lambda_i^2} +
 \frac{1- \lambda_j^2}{1+ \lambda_j^2}
= \frac{2(1- \lambda_i^2 \lambda_j^2) }{(1+ \lambda_i^2)(1+
\lambda_j^2)} \; .
\end{equation}
Therefore the area decreasing condition $|\lambda_i \lambda_j| <
1$ for $i\neq j$ is equivalent
 to the two-positivity of $S$, i.e. the sum of any two
eigenvalues is positive. We remark that curvature operator being
two-positive is preserved by the Ricci flow, see Chen \cite{ch} or
Hamilton \cite{ha3} for detail.

The two-positivity of a symmetric two tensor $P$ can be related to
the convexity of another tensor $P^{[2]}$ associated with $P$. The
following notation is adopted from Caffarelli-Nirenberg-Spruck
\cite{cns}. Let $P$ be a self-adjoint operator on an
$n$-dimensional inner product space. From $P$ we can construct a
new self-adjoint operator
\[
P^{[k]} = \sum_{i = 1}^k 1\otimes\cdots\otimes
\underset{i}{P}\otimes\cdots\otimes 1
\]
acting on the exterior powers $\Lambda^k$ by
\[
P^{[k]} (\omega_1 \wedge \cdots \wedge\omega_k) = \sum_{i = 1}^k
\omega_1 \wedge \cdots\wedge P(\omega_i)\wedge \cdots\wedge
\omega_k \; .
\]
With the definition of $P^{[k]}$, we have the following lemma.
\begin{lem}\label{lemS[2]}
Let $\mu_1
\le \mu_2 \le \cdots \le
\mu_n$ be the eigenvalues of  $P$
with corresponding
 eigenvectors $v_1 \cdots v_n$. Then
 $ P^{[k]} $ has eigenvalues $\mu_{i_1} +
\cdots +\mu_{i_k}$ and eigenvectors $v_{i_1}  \wedge \cdots\wedge
v_{i_k}$, $i_1 < i_2 \cdots < i_k$.
\end{lem}

Recall that the Riemannian metric $g$ and $S$ are both in
$T\Sigma\odot T\Sigma$, the space of symmetric two tensor on
$\Sigma$. We can identify  $S$ with a self-adjoint operator on the
tangent bundle through the metric $g$. Therefore $S^{[2]}$ and
$g^{[2]}$ are both sections of $(\Lambda^2(T\Sigma))^*\odot
\Lambda^2(T\Sigma)$ associated to $S$ and $g$ respectively. We
shall use orthonormal frames in the following calculation; this
has the advantage that $g$ is the identity matrix and we will not
distinguish lower index and upper index. With the above
interpretation and (\ref{eigs}), we have the following lemma.
\begin{lem}\label{2-positivity}
The area decreasing condition is equivalent to the convexity of
$S^{[2]}$.
\end{lem}
To show that the area decreasing condition is preserved, it
suffices to prove that the convexity of $S^{[2]}$ is preserved. In
fact, we prove the stronger result that
 the convexity of $S^{[2]}- \epsilon g^{[2]}$ for $\epsilon > 0$
 is preserved.

We  compute the evolution equation of  $S^{[2]}-\epsilon g^{[2]} $
in terms of the evolving orthonormal frames $\{F_a\}_{a=1\cdots
n}$ introduced earlier in (\ref{onf}). We will use indices $a, b,
\cdots $ to denote components in the evolving frames.
Denote $S_{ab} =  S(F_a,F_b )$ and $g_{ab}= g(F_a,F_b ) =
\delta_{ab}$.  Since $\{F_a \wedge F_b\}_{a<b}$ form a basis for $\Lambda^2
T\Sigma$, we have
\begin{equation}\label{S[2]}
\begin{split}
S^{[2]}(F_a \wedge F_b) = & \ S(F_a) \wedge F_b +  F_a \wedge
S(F_b)
= S_{ac} F_c \wedge F_b + F_a \wedge S_{ac}F_c \\
 =  & \ \sum_{c < d}(S_{ac} \delta_{bd} +  S_{bd} \delta_{ac}
 -S_{ad} \delta_{bc} -  S_{bc} \delta_{ad}) F_c \wedge F_d
 \, \quad {\rm and }\\
g^{[2]}(F_a \wedge F_b)=  & \ \sum_{c < d}( 2\delta_{ac}
\delta_{bd} - 2
\delta_{ad} \delta_{bc}) F_c \wedge F_d \, .
\end{split}
\end{equation}
We denote $S^{[2]}_{(ab) (cd)} = (S_{ac} \delta_{bd} +  S_{bd}
\delta_{ac}-S_{ad} \delta_{bc} -  S_{bc} \delta_{ad})$  and $g^{[2]}_{(ab)
(cd)} = 2\delta_{ac}
\delta_{bd} - 2
\delta_{ad} \delta_{bc}$. Thus the evolution equation of
$S^{[2]}-\epsilon g^{[2]} $ in terms of the evolving orthonormal
frames is
\begin{equation}\label{evol_two_p}
\begin{split}
 & \heat (S_{ac} \delta_{bd} +  S_{bd}
\delta_{ac}-S_{ad} \delta_{bc} -  S_{bc} \delta_{ad} - 2\epsilon \delta_{ac} \delta_{bd}
 + 2\epsilon  \delta_{ad} \delta_{bc} )
\\ & =  \
R_{e a e \alpha}S_{\alpha c}\delta_{bd} + R_{e c e
\alpha}S_{\alpha a}\delta_{bd} + R_{e b e \alpha}S_{\alpha
d}\delta_{ac} +
R_{e d e \alpha}S_{\alpha b}\delta_{ac} \\
& - R_{e a e \alpha}S_{\alpha d}\delta_{bc} - R_{e d e
\alpha}S_{\alpha a}\delta_{bc} - R_{e b e \alpha}S_{\alpha
c}\delta_{ad} - R_{e c e \alpha}S_{\alpha b}\delta_{ad}
\\
& + \  h_{\alpha e f} h_{\alpha e a}S_{fc} \delta_{bd} +
     h_{\alpha e f} h_{\alpha e c}S_{fa} \delta_{bd} +
     h_{\alpha e f} h_{\alpha e b}S_{fd} \delta_{ac} +
     h_{\alpha e f} h_{\alpha e d}S_{fb} \delta_{ac}  \\
     & - \  h_{\alpha e f} h_{\alpha e a}S_{fd} \delta_{bc} -
     h_{\alpha e f} h_{\alpha e d}S_{fa} \delta_{bc} -
     h_{\alpha e f} h_{\alpha e b}S_{fc} \delta_{ad} -
     h_{\alpha e f} h_{\alpha e c}S_{fb} \delta_{ad}  \\
& - \ 2h_{\alpha e a} h_{\beta e c}S_{\alpha \beta} \delta_{bd}
    -2h_{\alpha e b} h_{\beta e d}S_{\alpha \beta} \delta_{ac}
    + 2h_{\alpha e a} h_{\beta e d}S_{\alpha \beta} \delta_{bc}
    + 2h_{\alpha e b} h_{\beta e c}S_{\alpha \beta} \delta_{ad}  \, .
\end{split}
\end{equation}
Now, we are ready to prove that the area decreasing condition is
preserved along the mean curvature flow.
\begin{lem}\label{mainl} Under the assumption of Theorem A, with $S$
defined in (\ref{sym})and $S^{[2]}$ defined in (\ref{S[2]}),
suppose there exists an $\epsilon > 0$ such that
\begin{equation}\label{areadeceq}
  S^{[2]}- \epsilon g^{[2]} \ge 0
\end{equation}
holds on the initial graph. Then this is preserved along the mean
curvature flow.
\end{lem}

\begin{proof}.
Set \[M_{\eta} =  S^{[2]}- \epsilon g^{[2]} + \eta t  g^{[2]}\; .
\] Suppose the mean curvature flow exists on $[0, T)$. Consider any
 $T_1 < T$, it suffices to prove that $M_{\eta}  > 0$ on
$[0,T_1]$ for all $\eta < \frac{\epsilon}{2T_1}$. If not, there
will be a first time $0 <  t_0 \le T_1$ where $M_{\eta} = S^{[2]}-
\epsilon g^{[2]} +\eta t  g^{[2]}$ is nonnegative definite and
 has a null eigenvector $V=V^{ab}F_a \wedge F_b$  at some point $x_0 \in \Sigma_{t_0}$.
 We extend $V^{ab}$ to a parallel tensor in a neighborhood of $x_0$ along geodesic emanating out of
$x_0$, and defined $V^{ab}$ on $[0,T)$ independent of $t$. Define
 $f  = \sum_{ a< b, c < d} V^{ab}M_{\eta_{(ab)
(cd)}}V^{cd}$, then by (\ref{S[2]}), $f$ equals

\[
\begin{split}
& \sum_{ a< b, c < d} \left( S_{ac} g_{bd} +  S_{bd} g_{ac}-S_{ad}
g_{bc} -  S_{bc} g_{ad} + 2( \eta t- \epsilon) ( g_{ac} g_{bd}
 -  g_{ad} g_{bc})
\right)V^{ab}V^{cd}\; .
\end{split}\]
At $(x_0,t_0)$, we have  $f = 0$, $\nabla f = 0$ and $\heat f \leq
0$ where $\nabla$ denotes the covariant derivative and $\Delta$
denotes the Laplacian on $\Sigma_{t_0}$.

  We may assume that
 at  $(x_0,t_0)$ the orthonormal frames  $\{F_a\}$ is given by  $\{e_i\}$
  in (\ref{onf1_1}).
In the following, we use the orthonormal basis $\{e_i\}$ to write
down the condition $f=0$ and $\nabla f=0$ at $(x_0,t_0)$. The
basis $\{e_i\}$  diagonalizes $S$ with eigenvalues $
\{\lambda_i\}$ and   we order $
\{\lambda_i\}$ such that
\[
 \lambda_1^2 \ge
 \lambda_2^2 \ge \cdots \ge
 \lambda_n^2
\]
and
\begin{equation}\label{ineqs}
 S_{nn}=\frac{1- \lambda_n^2}{1+ \lambda_n^2}  \ge
 \cdots \ge
 S_{22}=\frac{1- \lambda_2^2}{1+ \lambda_2^2} \ge
 S_{11}=\frac{1- \lambda_1^2}{1+ \lambda_1^2}  \;.
 \end{equation}
 It follows from Lemma
(\ref{lemS[2]}) that $\{e_i\wedge e_j\}_{i < j}$ are the
eigenvectors of $M_{\eta}$. Thus we may assume that
\begin{equation}
V=e_1\wedge e_2 .
\end{equation}

 At $(x_0,t_0)$, the condition $f = 0$ is the same as
 \begin{equation}\label{S1122}
 S_{11} + S_{22} =  2\epsilon - 2 \eta t_0 > 0\; .
 \end{equation}
This is equivalent to
\[\frac{2(1- \lambda_1^2 \lambda_2^2) }{(1+ \lambda_1^2)(1+
\lambda_2^2)} = 2(\epsilon - \eta t_0) > 0.\]

Thus

\begin{equation}\label{f=0} \lambda_1\lambda_2 <1 \,\,\text{and}\,\, \lambda_i<1
\,\,\text{for}\,\, i\ge 3\; .\end{equation}
 Next, we
 compute the covariant derivative of the restriction of $S$
on $\Sigma$.
\[
\begin{split}
& (\nabla_{e_i}S)(e_j,e_k)\\
= \quad & e_i(S(e_j,e_k))- S(\nabla_{e_i}e_j,e_k)-S(e_j,\nabla_{e_i}e_k) \\
= \quad & S(\nabla_{e_i}^{M}e_j -\nabla_{e_i}e_j,e_k)+
S(e_j,\nabla_{e_i}^{M}e_k-\nabla_{e_i}e_k)\\
= \quad & h_{\alpha i j}S_{\alpha k} + h_{\beta i k}S_{\beta j}\; .
\end{split}
\]
So
\[
S_{jk,i }=h_{\alpha ij }S_{\alpha k} +
 h_{\beta ik}S_{\beta j} \; .
\]
Recall that  $V_{ab}$ is  parallel at $(x_0,t_0)$ , $V^{12}
 =1$ and all other components of $V^{ab}$ is zero.
 At $(x_0,t_0)$, $\nabla f = 0$ is equivalent to
\[
\begin{split}
 \ 0 = & \ \ \sum_{ i< j, k < l}  \nabla_{e_p} ((S_{ik} \delta_{jl} +  S_{jl} \delta_{ik}
 -S_{il} \delta_{jk} - S_{jk} \delta_{il}+2(\eta t - \epsilon) (\delta_{ik}
\delta_{jl}
- \epsilon  \delta_{il}\delta_{jk}))V^{ij}V^{kl}) \\
= & \ \nabla_{e_p} S_{11}+ \nabla_{e_p} S_{22}\\
= & \ 2h_{\alpha p 1}S_{\alpha 1} +
 2 h_{\beta p 2}S_{\beta 2}\; .
\end{split}
\]
Since $S_{n+q, l} =
 - \frac{2\lambda_q \delta_{ql}}{1+\lambda_q^2}$, we have
\begin{equation}\label{gradient}
\frac{\lambda_1}{1+\lambda_1^2} h_{n+1, p 1} +
\frac{\lambda_2}{1+\lambda_2^2} h_{n+2, p 2} = 0
\end{equation}
for any $p$.

By (\ref{evol_two_p}), at $(x_0,t_0)$,  we have
\begin{equation}\label{null_c2}
 \begin{split}
   \heat f
 & =  \ 2 \eta + \  2 R_{k 1 k \alpha}S_{\alpha 1} +
 2 R_{k 2 k \alpha}S_{\alpha 2}
    + 2h_{\alpha k j} h_{\alpha k 1}S_{j1}+
      2h_{\alpha k j} h_{\alpha k 2}S_{j2} \\
 & - \ 2h_{\alpha k 1} h_{\beta k 1}S_{\alpha \beta}
  -2 h_{\alpha k 2} h_{\beta k 2}S_{\alpha \beta}. \\
 \end{split}
\end{equation}
 The ambient curvature term can be calculated using Lemma 4.2 and we
 derive
\begin{equation}\label{null_c3}
\begin{split}
&\sum_{k,\alpha}  R_{k 1 k \alpha}S_{\alpha 1} +
 R_{k 2 k \alpha}S_{\alpha 2} . \\
 = & \ {(k_1-k_2)(n-1)}
 \sum_{i=1}^2\frac{\lambda_i^2}{(1+\lambda_i^2)^2} +(k_1 + k_2)
 \sum_{i=1}^2\frac{\lambda_i^2}{(1+\lambda_i^2)^2}
\left[\sum_{j\not= i}\frac{1-\lambda_j^2}{(1+\lambda_j^2)}
\right]\; .
\end{split}
\end{equation}

This can be simplified as
 \begin{equation}
\begin{split}
  & \ {(k_1-k_2)(n-1)}
 \sum_{i=1}^2\frac{\lambda_i^2}{(1+\lambda_i^2)^2} +(k_1 + k_2)
 \sum_{i=1}^2\frac{\lambda_i^2}{(1+\lambda_i^2)^2}
\left[\sum_{j > 3}\frac{1-\lambda_j^2}{(1+\lambda_j^2)} \right] \\
+
 & \ (k_1 + k_2)\left[
\frac{\lambda_1^2}{(1+\lambda_1^2)^2}\frac{1-\lambda_2^2}{(1+\lambda_2^2)}
 + \frac{\lambda_2^2}
{(1+\lambda_2^2)^2}\frac{1-\lambda_1^2}{(1+\lambda_1^2)}\right]\\
 = & \ {(k_1-k_2)(n-1)}
 \sum_{i=1}^2\frac{\lambda_i^2}{(1+\lambda_i^2)^2} +(k_1 + k_2)
 \sum_{i=1}^2\frac{\lambda_i^2}{(1+\lambda_i^2)^2}
\left[\sum_{j > 3}\frac{1-\lambda_j^2}{(1+\lambda_j^2)} \right] \\
+
 & \ (k_1 + k_2)\left[
 \frac{(\lambda_1^2+ \lambda_2^2) (1- \lambda_1^2\lambda_2^2)}
{(1+\lambda_1^2)^2(1+\lambda_2^2)^2}\right] \;.
\end{split}
\end{equation}
This is nonnegative by equation (\ref{f=0}).

Using the relations in (\ref{block}) again, the last four terms on
the right hand side of (\ref{null_c2}) can be rewritten as
\begin{equation}\label{null_c4}
\begin{split}
\ & \sum_{p, k}2h_{n+p, k 1}^2 S_{11}+
      2h_{n+p, k 2}^2 S_{22}
  +2 h_{n+p, k 1}^2 S_{pp}
  +2 h_{n+p, k 2}^2 S_{pp} \\
& = \ \sum_{k} (2h_{n+1, k 1}^2S_{1 1} + 2h_{n+2, k 1}^2S_{1 1}+
 2h_{n+1, k 2}^2S_{2 2} + 2h_{n+2, k 2}^2S_{2 2} \\
 & + \  2h_{n+1, k 1}^2S_{11} +2h_{n+2, k 1}^2S_{22} +
  2h_{n+1, k 2}^2S_{11} + 2h_{n+2, k 2}^2S_{22} ) \\
& + \ \sum_{q  \ge 3,k}2h_{n+q, k 1}^2S_{1 1}+ 2h_{n+q, k 2}^2S_{2
2} + 2h_{n+q, k 1}^2S_{qq} + 2h_{n+q, k 2}^2S_{qq}\; .
\end{split}
\end{equation}

 Since $S_{ii}
\ge S_{11}$ for $i \ge 2$, it is clear that (\ref{null_c4})
 is  nonnegative if  $S_{11} \ge 0 $.
Otherwise, from (\ref{S1122}), we may assume that
\begin{equation}\label{eqs}
S_{11} < 0,\  S_{22} > 0 \ {\rm and } \ S_{11} + S_{22}  > 0 \; .
\end{equation} In particular,
we have $\lambda_2^2 < \lambda_1^2$ and $\lambda_1^2 \lambda_2^2 <
1$. From (\ref{gradient}), we have

\[  h_{n+1 ,p 1}^2
=\frac{\lambda_2^2(1+\lambda_1^2)^2}{\lambda_1^2(1+\lambda_2^2)^2}
h_{n+2, p 2}^2 .\]

Since  $\lambda_2^2 < \lambda_1^2$ and  $\lambda_1^2 \lambda_2^2 <
1$, we have
$\frac{\lambda_2^2(1+\lambda_1^2)^2}{\lambda_1^2(1+\lambda_2^2)^2}
< 1$. Thus
\begin{equation}\label{secineq}
 h_{n+1, p 1}^2 \le   h_{n+2, p 2}^2
{\rm \ for \ all \ } p \ \ge 1.\end{equation} Recall that $S_{qq}
\ge S_{22}$ for $q \ge 3$. The right hand side of (\ref{null_c4})
can be regrouped as
\[
\begin{split} & \ \sum_{k} \Big[ ( 4 h_{n+1, k 1}^2S_{1 1}+ 4h_{n+2, k 2}^2S_{22})
 +  2h_{n+2, k 1}^2(S_{1 1}+S_{22}) +
 2h_{n+1, k 2}^2(S_{11} +S_{2 2}) \Big] \\
& + \ \sum_{q  \ge 3,k}\Big[2h_{n+q, k 1}^2(S_{1 1}+S_{qq}) +
2h_{n+q, k
2}^2( S_{2 2} + S_{qq})\Big] \; . \\
\end{split}
\]
This  is  nonnegative by (\ref{ineqs}) ,(\ref{eqs}), and
(\ref{secineq}). Thus, we have $\heat f \ge 2\eta
> 0$ at $(x_0,t_0)$ and this is a contradiction.
\end{proof}

\noindent Remark: The condition $ S^{[2]}- \epsilon g^{[2]} \ge 0$
is equivalent to  $\frac{(1- \lambda_i^2 \lambda_j^2) }{(1+
\lambda_i^2)(1+ \lambda_j^2)} \ge  \epsilon$ for all $i \neq j$.
In particular, we have $\lambda_i^2 \le
\frac{1-\epsilon}{\epsilon}$. This implies that the Lipschitz norm
of $f$ is preserved along the mean curvature flow.

\section{Long time existence and convergence}
In this section, we prove Theorem A using the evolution equation
(\ref{lneq}) of $\ln *\Omega$.

\begin{proof}{\,\,\it of Theorem A}.
Since $|\lambda_{i}\lambda_{j}| < 1$  for $ i \neq j$ and
$\Sigma_1$ is compact, we can find an $\epsilon > 0$ such that
$\frac{(1- \lambda_i^2 \lambda_j^2) }{(1+ \lambda_i^2)(1+
\lambda_j^2)} \ge  \epsilon$ for all $i \neq j$. By Lemma
(\ref{mainl}), the condition $\frac{(1- \lambda_i^2 \lambda_j^2)
}{(1+ \lambda_i^2)(1+ \lambda_j^2)} \ge  \epsilon$ for all $i \neq
j$ is preserved along the mean curvature flow. In particular, we
have $|\lambda_{i}\lambda_{j}| \le \sqrt{1-\epsilon}$ and
$\lambda_i^2\le\frac{1-\epsilon}{\epsilon}$. This implies
$\Sigma_t$ remains the graph of a map $f_t:\Sigma_1\rightarrow
\Sigma_2$ whenever the flow exists. Each $f_t$ has uniformly
bounded $|d f_t|$.

We look at the evolution equation (\ref{lneq}) of $\ln *\Omega$.
The quadratic terms of the second fundamental form in equation
(\ref{lneq}) is
\[
\begin{split}
& \sum_{\alpha, i,k} h_{\alpha ik}^2 + \sum_{k, i}\lambda_{ i}^2
h_{n+i,ik}^2
+2\sum_{k, i<j}\lambda_{i}\lambda_{j} h_{n+j, ik} h_{n+i, jk}\\
 = & \ \delta |A|^2  +
\sum_{k, i}\lambda_{ i}^2 h_{n+i,ik}^2 +(1- \delta) |A|^2
+2\sum_{k, i<j}\lambda_{i}\lambda_{j} h_{n+j, ik} h_{n+i, jk}\; .
\end{split}
\]
 Let $1-\delta = \sqrt{1-\epsilon}$. Using $|\lambda_i\lambda_j|\leq 1-\delta$, we conclude that this
term is  bounded below by $\delta |A|^2$ .

 By equation (\ref{useful_identity}), the curvature term in (\ref{lneq}) equals
\begin{equation}\label{positive1}
 \begin{split}
 & \ \frac{(k_1-k_2)(n-1)}{2}
 \sum_{i=1}^n \frac{\lambda_{ i}^2}{1+\lambda_i^2} +
(k_1 +k_2) \sum_{i=1}^n \frac{\lambda_{ i}^2}{1+\lambda_i^2}
\left[\sum_{j\not=
i}\frac{1-\lambda_j^2}{2(1+\lambda_j^2)}\right]\; .
 \end{split}
\end{equation}

 The second term on the
right hand side of (\ref{positive1}) can be simplified as
\begin{equation}\label{positive2}
\begin{split}
 & \ \sum_{i=1}^n\frac{\lambda_{ i}^2}{1+\lambda_i^2}
\left[\sum_{j\not= i}\frac{1-\lambda_j^2}{2(1+\lambda_j^2)}
 \right]
=  \  \sum_{i=1}^n\sum_{ i \neq j} \frac{\lambda_{ i}^2  -
\lambda_{ i}^2 \lambda_{ j}^2} {2(1+\lambda_i^2)(1+\lambda_j^2)}
 \\
= & \ \sum_{ i < j} \frac{\lambda_{ i}^2 + \lambda_{ j}^2 - 2
\lambda_{ i}^2 \lambda_{ j}^2} {2(1+\lambda_i^2)(1+\lambda_j^2)}.
\end{split}
\end{equation}

This is non-negative because $|\lambda_{i}\lambda_{ j}| \le 1-
\delta$. Thus $\ln*\Omega$ satisfies the following differential
inequality with $k_1 \ge |k_2|$:
\begin{equation}\label{lneq2}
\frac{d}{dt}\ln *\Omega\geq \Delta \ln *\Omega +\delta |A|^2\; .
\end{equation}
According to the maximum principle for parabolic equations,
$\min_{\Sigma_t} \ln *\Omega $ is nondecreasing in time. In
particular, $*\Omega \ge min_{\Sigma_0}*\Omega= \Omega_0$ is
preserved and $*\Omega$ has a positive lower bound.  Let $u =
\frac{\ln *\Omega - \ln \Omega_0 + c} {- \ln \Omega_0+c}$ where
$c$ is a positive number such that $- \ln \Omega_0+c > 0$. Recall
that $ 0 < *\Omega \le 1$. This implies that $0 < u \le 1$ and $u$
satisfies the following differential inequality
\[\frac{d}{dt}u\geq \Delta u +\frac{\delta}{- \ln \Omega_0 +c}
|A|^2.
\]

Because $u$ is also invariant under parabolic dilation, it follows
from the blow-up analysis in  the proof of Theorem A \cite{wa3}
that the mean curvature flow of the graph of $f$
 remains a graph and exists for all time under the assumption that $k_1 \ge |k_2|$.

Using $\lambda_i^2\le\frac{1-\epsilon}{\epsilon}$ and
 $\lambda_i\lambda_j\le \sqrt{1-\epsilon}$, it is not hard to show

\begin{equation}\label{curineq}
 (k_1+k_2)\sum_{ i < j}
\frac{\lambda_{ i}^2 + \lambda_{ j}^2 - 2 \lambda_{ i}^2 \lambda_{
j}^2} {2(1+\lambda_i^2)(1+\lambda_j^2)} \ge c_1\sum_{i=1}^n
{\lambda_i^2} \geq c_1\ln\prod_{i=1}^n(1+\lambda_i^2)\;
\end{equation}
where $c_1$ is a constant that depends on $\epsilon, k_1$ and
$k_2$.

Recall equation (\ref{Omega}) and we obtain
\[\frac{d}{dt} \ln *\Omega\geq \Delta
\ln *\Omega-c_3\ln*\Omega\; .\]

By the comparison theorem for parabolic equations,
$\min_{\Sigma_t} \ln *\Omega$ is non-decreasing in $t$ and
 $\min_{\Sigma_t} \ln*\Omega \rightarrow 0$ as $t\rightarrow \infty$.
 This implies that $\min_{\Sigma_t} *\Omega \rightarrow 1$ and $\max |\lambda_i|\rightarrow 0$ as $t\rightarrow
 \infty$. We can then apply Theorem B in  \cite{wa3} to conclude smooth convergence to a constant map at infinity.

\end{proof}

\end{document}